\documentclass[CJK,11pt,twoside, openany]{ctexart}
\usepackage{amsfonts}
\usepackage{}
\usepackage{latexsym}
\usepackage{amsthm}
\usepackage{bbm}
 \usepackage{mathrsfs}
 \usepackage{fancyhdr}
\usepackage{mathrsfs}
\usepackage{cancel}
\usepackage{color}
\usepackage{multirow}
\usepackage{eso-pic}
\usepackage{dsfont}
\usepackage{amssymb}
\usepackage{CJKspace}

\def\thebibliography#1{\center{\bf\normalsize References}\list
 {[\arabic{enumi}]}{\settowidth\labelwidth{[#1]}\leftmargin\labelwidth
 \advance\leftmargin\labelsep
 \usecounter{enumi}}
 \def\newblock{\hskip .11em plus .33em minus .07em}
 \sloppy\clubpenalty4000\widowpenalty4000
 \sfcode`\.=1000\relax}

\makeatletter
\def\cleardoublepage{\clearpage\if@twoside \ifodd\c@page\else
   \hbox{}\thispagestyle{empty}\newpage\addtocounter{page}{-1}
   \if@twocolumn\hbox{}\newpage\fi\fi\fi}
\makeatother
\usepackage{graphicx}
\def\nnb{\nonumber}
\def\ds{\displaystyle}
\def\cd{\cdot}

\def\all{  \, \forall \, }

\newcommand{\refeq}[1]{~$(\ref{#1})$}
\newcommand{\eqref}[1]{~$(\ref{#1})$ }
\newcommand{\thb}[1]{~{\rm (#1)}}
\def\endpf{\hfill$\Box$\vspace{0.4cm}}

\def\eqae{ \, {\rm a.e. } \,\, }

\def\Gl{\lambda}

\def\Go{\omega}

\def\Ga{\alpha}
\def\Gb{\beta}

\def\Gd{\delta}

\def\Gve{\varepsilon}

\def\Gt{\theta}
\def\Gvp{\varphi}

\def\GT{\Theta}

\def\GO{\Omega}
\def\bt{{\bar t}}
\def\bu{{\bar u}}

\def\by{{\bar y}}

\def\mcP{{\mathscr P}}

\def\mcU{{\mathscr U}}

\def\bpsi{{\bar \psi}}

\def\hu{\hat u}

\def\hx{\hat x}
\def\hy{\hat y}

\def\hGt{\hat \Gt}

\def\tis{\tilde s}

\def\tiu{\tilde u}

\def\tix{\tilde x}
\def\tiy{\tilde y}

\def\tiGO{\tilde \GO}

\def\tiGt{\tilde \Gt}

\def\qq{\qquad}
\def\q{\quad}
\newcommand{\dpp}[2]{{\partial {#1} \over \partial {#2}}}
\newcommand{\doo}[2]{{d {#1} \over d {#2}}}

\newcommand{\pri}{{\prime}}
\newcommand{\prii}{{\prime\prime}}

\def\pa{\partial}

\def\lra{\longrightarrow}

\def\limsup{\mathop{\overline{\rm lim}}}
\def\liminf{\mathop{\underline{\rm lim}}}
\def\essinf{\mathop{{\rm essinf}}}
\def\esssup{\mathop{{\rm esssup}}}

\def\IR{\mathds{R}}

\newcommand{\set}[1]{\left\{#1\right\}}
\newcommand{\ip}[1]{\left\langle #1\right\rangle}

\newtheorem{Definition}{Definition}[section]
\newtheorem{Theorem}[Definition]{Theorem}
\newtheorem{Lemma}[Definition]{Lemma}
\newtheorem{Corollary}[Definition]{Corollary}
\newtheorem{Proposition}[Definition]{Proposition}
\newtheorem{Remark}{Remark}[section]
\newtheorem{Example}{Example}

\begin{document}
\title{\bf  Optimal Blowup Time for Controlled Ordinary Differential Equations \thanks{This work was supported in part by 973 Program (No. 2011CB808002)  and  NSFC (No. 11371104).}}

\author{Hongwei Lou\footnote{School
of Mathematical Sciences, and LMNS, Fudan University, Shanghai
200433, China (Email: \texttt{hwlou@fudan.edu.cn}).}
~~and~~Weihan Wang\footnote{School of Mathematical Sciences, Fudan University, Shanghai 200433, China (Email: \texttt{11210180039@fudan.edu.cn}).}
}

\date{}

\maketitle

\begin{quote}
\footnotesize {\bf Abstract.}  Both the shortest and the longest blowup time for a controlled system are
considered. Existence result and maximum principle for optimal triple are established.
Thanks to some monotonicity of the controlled system, some kinds of ``the front part local optimality" for optimal triple is established. Then proofs of the main results become easy, clear and abundant.

\textbf{Key words and phrases.}  Optimal blowup time,  Front part local optimality, Existence, Maximum principle

\textbf{AMS subject classifications.} 49J15, 34A34
\end{quote}

\normalsize

\def\theequation{1.\arabic{equation}}
\setcounter{equation}{0} 
\setcounter{Definition}{0} \setcounter{Remark}{0}
 \vspace{6mm}


 \section{Introduction.}

The blowup phenomenon of evolution equations is of great importance, and has been studied by numerous researchers. A typical model (c.f. \cite{Bandle}) is described as $u_t-\Delta u=f(x,t,u,\nabla u)$, which characterizes the temperature of a substance in a chemical reaction. The blowup phenomenon can represent a dramatic increase in temperature which leads to the ignition of a chemical reaction.
In the past 50 years, most of the former works focus on the existence of blowup solutions and the blowup rates (see \cite{Escobedo}, \cite{Glassey}, \cite{Guo}, \cite{Yordanov} and \cite{Zhang}, for examples). It is natural to ask what is a good/best method to control the blowup time. In certain cases, we hope to minimize the blowup time, while we may want to maximize the blowup time in other cases. So it is meaningful to consider the relevant time optimal control problems. As Barron and Liu mentioned in their paper(see \cite{Barron}), although the researchers' initial interest is about the optimal control to the distributed systems, they met some difficulties. Hence, some researchers discuss the relevant problems governed by ordinary differential equations.


According to our knowledge, the related research is very limited. Barron and Liu \cite{Barron} posed an optimal control problem to maximize the blowup time in 1996 . They consider an autonomous system, which is described as:
\begin{equation}\label{E101}
   \left\{\begin{array}{ll}\ds  \doo {y(t)} t=f(y(t),u(t)),&　t>0,\\
   y(0)=x. & \end{array}\right.
\end{equation}
where for some $p>1$,
\begin{equation}\label{E102}
   \frac{x \cd f(x,z)}{|x|^{p+1}}\to  1,\qq  \mbox{as }\q |x|\to \infty
\end{equation}
uniformly in $z$.

For a fixed control $u(\cdot )$, the blowup time is considered as a map $x\longmapsto T_x(u)$. And the value function is defined as:
\begin{equation}\label{E103}
   V(x)=\sup_u T_x(u),\qq V:  \IR^n\lra [0,\infty].
\end{equation}
They study the properties of the blowup time and the value function. Then they conclude from the dynamic programming principle that the value function is the unique continuous viscosity solution of the Hamilton-Jacobi equation:
\begin{equation}\label{E104}
   1+\max_z D_xV(x)\cdot f(x,z)=0
\end{equation}
which leads to the maximum principle.

In 2011,  Lin and  Wang \cite{LinWang} studied an optimal control problem to minimize the blowup time, which is governed by a special non-autonomous system:
\begin{equation}\label{E105}
   \left\{\begin{array}{ll}\ds  \doo {y(t)} t=|y(t)|^{p-1}y(t)+B(t)u(t),&　t>0,\\
   y(0)=y_{0}. & \end{array}\right.
\end{equation}
The main results are the existence and the maximum principle of the optimal control problem. Their strategy to get the existence is to transform the problem into the classical case. In detail, they show the existence of a series of relevant problems $(P_{R})$, where the target sets are the sphere of the ball in $\IR^{n}$, centered at the origin and of different radius $R$. And prove the existence of the origin problem by taking limit of $R$. Then, they introduce a new penalty function to conclude the maximum principle. In \cite{LWX}, Lou, Wen and Xu gave another approach to discuss the problem.

In this paper, we will consider the minimal/maximal blowup time optimal control problem, governed by a general system, which covers the systems mentioned above:
\begin{equation}\label{E106}
   \left\{\begin{array}{ll}\ds  \doo {y(t)} t=f(t,y(t),u(t)),&　t>0,\\
   y(0)=y_0. & \end{array}\right.
\end{equation}
where  $f(t,y,u)$ takes the following form:
\begin{equation}\label{E107}
f(t,y,u)=G(t,|y|){y\over |y|}+A(t)y +b(t,u), \qq\all (t,y,u)\in [0,+\infty)\times \IR^n\times U
\end{equation}
with  $(U,\rho)$ being a separable metric  space, $G(\cd,\cd)$ being a function on $[0,+\infty)^2$, $A(\cd)$ being  an $n\times n$-matrix-valued function on  $[0,+\infty)$ and $b(\cd,\cd)$ being an
$n$ dimensional vector-valued  function on $[0,+\infty)\times U$.
\if{
We say $y(\cd)\in C([0,T);\IR^n)$ is a solution of \refeq{E106} on $[0,T)$ if $y(\cd)$ satisfies
$$
y(t)=y_0+\int^t_0f(s,y(s),u(s))\, ds, \qq t\in [0,T).
$$
The time optimal control problem considered can be described as
$$
\textbf{Problem (TI)}:  \qq \inf \left\{T\in (0,+\infty)\left|\begin{array}{l}\ds   y(\cd;u(\cd)) \, \mbox{exists on \,} [0,T),\\
\ds \lim_{t\to T^-}|y(t;u(\cd))|=+\infty, \end{array}\q  u(\cdot )\in \mcU_{ad}\right. \right\},
$$
其中状态函数 $y(\cd)=y(\cd; u(\cd))$ 表示方程 \refeq{E106} 对应于控制 $u(\cd)$ 的解, 状态函数 $y(\cd)$ 是取值于 $\IR^n$ 的向量值函数, 而控制函数 $u(\cd)$ 取值于 $\IR^m$ 的非空子集 $U$, $\mcU_{ad}$  为 $\ds \mcU\equiv \set{u(\cd): [t_0,+\infty )\to U| u(\cd)\, \mbox{measurable}}$ 的满足一定条件的子集.
}\fi
We say \refeq{E106} holds on $[0,T)$ or $y(\cd)$ is a solution of \refeq{E106} on $[0,T)$ for some $T>0$  always means that $y(\cd)\in C[0,T)$
and
$$
y(t)=y_0+\int^t_0f(s,y(s),u(s))\, ds, \qq\all t\in (0,T).
$$
Denote
\begin{equation}\label{E108}
\begin{array}{l}\ds
\mcU=\set{u(\cd):[0,+\infty)\to U\Big|\, u(\cd)\q\mbox{is measurable}\,},\\
\mcP=\set{(T,y(\cd),u(\cd))\in (0,+\infty)\times C([0,T);\IR^n)\times \mcU\Big| \mbox{\refeq{E106} \,holds on }\, [0,T)},\\
\ds \mcP_{ad}=\set{(T,y(\cd),u(\cd))\in \mcP\Big| \lim_{t\to T^-}|y(t)|=+\infty},\\
\ds \mcU_{ad}=\set{u(\cd)\in \mcU\Big| (T,y(\cd),u(\cd))\in \mcP_{ad}}.
\end{array}
\end{equation}
Moreover, $\mcP$, $\mcP_{ad}$ and $\mcU_{ad}$ are named as the set of feasible triples, the set of admissible triples and the set of admissible controls, respectively.

If $\mcU_{ad}\ne \emptyset$, the corresponding minimal time optimal control problem is:

\textbf{Problem (TI)}: Find $(\bt,\by(\cd),\bu(\cd))\in \mcP_{ad}$  such that
\begin{equation}\label{E109}
\bt=\inf_{(T,y(\cd),u(\cd))\in \mcP_{ad}}T.
\end{equation}

If $\mcU_{ad}=\mcU$, then we can consider the maximal time optimal control problem:

\textbf{Problem (TS)}: Find $(t^*,y^*(\cd),u^*(\cd))\in \mcP_{ad}$  such that
\begin{equation}\label{E110}
t^*=\sup_{(T,y(\cd),u(\cd))\in \mcP_{ad}}T.
\end{equation}

The focus of this paper is to introduce a new approach to yield the maximum principle. The key of our strategy is to bridge the gap between classical cases and the blowup ones. Roughly speaking,  our results will be based on establishing ``the front part local optimality" by some monotonicity of controlled systems.
More precisely, if $\bu(\cd)$ is an optimal control for a time optimal control problem and $\bt$ is the optimal time, then its rear part is also optimal, i.e., for any $T\in (0,\bt\,)$,  $\bu(\cd)|_{[T,\bt\,)}$ is also an optimal control for the  time optimal control problem restricted on $[T,\bt\,]$. However, for a non-autonomous system, it is not necessary that $\bu(\cd)|_{[0,T]}$ is an optimal control for the  time optimal control problem in the front part.
Nevertheless, by studying some monotonicity of the controlled system, we  can construct some kind of  the front part local optimality
of the optimal trajectory before blowup. Then the maximum principle follows by taking limit of the classical results. On the other hand, as to autonomous systems, translation invariance of the trajectory ensures the local optimality---``the front part local optimality" as well as ``the rear part local optimality". Thus, results for autonomous systems can be got much easily and under relatively weaker assumptions than non-autonomous cases.

Based on the new approach,  the controlled systems we considered are more general than those considered in the previous works.

The existence results (see Theorems \ref{T202a} and \ref{T203}) will be established in Section 2. Section 3 will be devoted to the maximum principles for optimal control to Problem (TI) (see Theorems \ref{T302} and \ref{T303}). While Section 4 is devoted to the maximum principles for optimal control to Problem (TS) (see Theorem \ref{T304}). In Section 5, we mention that for autonomous systems, the maximum principles for optimal triples are relatively easy to be established. Finally, we will list some examples to show that our results can be applied to most of interested systems.
\bigskip

\bigskip

\def\theequation{2.\arabic{equation}}
\setcounter{equation}{0} 
\setcounter{Definition}{0} \setcounter{Remark}{0}\section{Existence of Time Optimal Control Problem }

In this section, we will discuss the existence of optimal control. We make the following assumptions:

(P1) Let $(U,\rho)$ being a separable metric  space;

(P2) Function $G(t,r)$ is measurable in $t\in [0,+\infty)$, continuously differentiable in $r\in [0,+\infty)$ and
\begin{equation}\label{E201}
G(t,0)=0,   \qq  \all  t \in [0,+\infty).
\end{equation}
Moreover,for any $M>0$,
\begin{equation}\label{E202}
\sup_{(t,r)\in [0,M]^2}|{\pa G(t,r)\over \pa r}|<+\infty.
\end{equation}

(P3) Let $A(\cd)\in L^\infty_{loc}([0,+\infty);\IR^{n\times n})$, i.e.,  for any $T>0$,
\begin{equation}\label{E203}
 \sup_{t\in [0,T]} \|A(t)\|<+\infty,
\end{equation}
where $\|A\|$ represents the norm of a $n\times n$ matrix $A$: $\|A\|=\sup_{x\in S^{n-1}}|Ax|$.

(P4) Function $b(\cd,\cd)$ takes values in $\IR^n$ and it is a Carath\'eodory function that it is measurable in the first variable and continuous in the second variable. Moreover, $U(t)\equiv \set{b(t,u)|u\in U}$ is a convex compact set.

(P5) There exists an $R_0>0$ and a nonnegative function $\zeta(\cd)$ defined on $[R_0,+\infty)$, which satisfies
\begin{equation}\label{E204}
G(t,r)-\|A(t)\|\,r-\sup_{u\in U} |b(t,u)|\geq \zeta(r),\q\all (t,r,u)\in [0,+\infty)\times [R_0,+\infty)\times U,
\end{equation}
\begin{equation}\label{E205}
\int^{+\infty}_{R_0}{1\over \zeta(r)}\, dr<+\infty.
\end{equation}

We notice that the solution of equation \refeq{E106} is well-posed before blow-up according to (P2). So it is reasonable to represent the solution of  \refeq{E106} as  $y(\cd;u(\cd))$.

We introduce the following lemma to establish the existence theorem:

\begin{Lemma}\label{T201}
Assume \thb{P1}---\thb{P4} hold. Let $(T,y(\cd),u(\cd))\in \mcP$, $\ds\limsup_{t\to T^-}|y(t)|<+\infty$, $u_k(\cd)\in \mcU$ and
\begin{equation}\label{E206}
b(\cd,u_k(\cd))\to b(\cd,u(\cd)), \qq \mbox{weakly in }\, L^2([0,T+1];\IR^n).
\end{equation}
Then, there exist $\Gd>0$ and $K>0$, such that $y_k(\cd)\equiv y(\cd;u_k(\cd))$, which is the solution of equation \refeq{E106} with control $u_k(\cd)$, exists on $[0,T+\Gd]$ and satisfies
\begin{equation}\label{E207}
|y_k(t)-y(t)|\leq 1, \qq\all t\in [0,T+\Gd]
\end{equation}
when $k\geq K$.
\end{Lemma}
\proof Let
$$
v_k(t)=b(t,u_k(t)), \q v(t)=b(t,u(t)), \qq t\in [0,T+1].
$$
By (P1)---(P4), $\ds\limsup_{t\to T^-}|y(t)|<+\infty$, and the basic theory of ordinary differential equation, there exists $\Gd\in (0,1)$, such that $y(\cd;u(\cd))$ exists on $[0,T+\Gd]$.

Let
$$
M_b=\max_{t\in [0,T+\Gd]\atop u\in U}|b(t,u)|, \q M_A=\max_{t\in [0,T+1]}\|A(t)|,
$$
$$
M=\max_{t\in [0,T+\Gd]} |y(t)|+T+\Gd+1, \q N=\sup_{(t,r)\in [0,M]^2}|{\pa G(t,r)\over \pa r}|.
$$
Let $\ds \Ga= {1\over 3 e^{(3N+M_A)(T+1)}}$, $\ell>{2(T+1)M_b\over  \Ga}$ and $\ell$ be an integer.
By the weak convergence, it is not difficult to prove that for some $K>0$,
\begin{equation}\label{E208}
\Big|\int^{j(T+1) \over \ell}_0 (v_k(t)-v(t))\, dt \Big| \leq \Ga,\q \all j=1,2,\ldots,\ell-1
\end{equation}
when $k\geq K$.
Thus,
\begin{equation}\label{E209}
\Big|\int^t_0 (v_k(t)-v(t))\, dt \Big| \leq 2\Ga,\q \all t\in [0,T+1], \q k\geq K.
\end{equation}
We claim that when $k\geq K$,
\begin{equation}\label{E210}
|y_k(t)|<M, \qq\all t\in [0,T+\Gd].
\end{equation}
Otherwise,  for some $k\geq K$, there exists $S\in (0,T+\Gd]$, such that $|y_k(S)|=M$,
\begin{equation}\label{E211}
|y_k(t)|<M, \qq\all t\in [0,S).
\end{equation}
We have
\begin{eqnarray*}
&& |y_k(t)-y(t)|\\
&=& \Big|\int^t_0\Big(G(s,|y_k(s)|){y_k(s)\over |y_k(s)|}-G(s,|y(s)|){y(s)\over |y(s)|}+A(s)(y_k(s)-y(s))+v_k(s)-v(s)\Big)\, ds  \Big|\\
&\leq &  \Big|\int^t_0\Big(G(s,|y_k(s)|)-G(s,|y(s)|)\Big){y_k(s)\over |y_k(s)|} \, ds  \Big|\\
& & + \Big|\int^t_0 {G(s,|y(s)|)\over |y(s)|} \Big(y_k(s)-y(s)+(|y(s)|-|y_k(s)|){y_k(s)\over |y_k(s)|} \Big)\, ds  \Big|\\
&& + M_A\int^t_0 |y_k(s)-y(s)|\, ds +2\Ga\\
&\leq &  (3N+M_A) \int^t_0 |y_k(s) -y(s)|  \, ds   +2\Ga,\qq\all t\in  [0,S].
\end{eqnarray*}
Adopting Grownwall's inequality, we can get
\begin{equation}\label{E212}
|y_k(t)-y(t)|\leq  3\Ga e^{(3N+M_A)t}\leq  1, \qq\all t\in [0,S].
\end{equation}
In particular,
$$
|y_k(S)|\leq |y(S)|+1<M,
$$
which contradicts $|y_k(S)|=M$. Therefore, \refeq{E210} holds. Further, we get \refeq{E204A} from \refeq{E210} (see \refeq{E212}).

The proof is completed.
\endpf

\begin{Lemma}\label{T202}
Assume \thb{P1}---\thb{P4} hold. Let $(T,y_k(\cd),u_k(\cd))\in \mcP$ satisfy
\begin{equation}\label{E213}
\lim_{k\to +\infty} |y_k(T)|=+\infty.
\end{equation}
Then, $\bt\leq T$, where $\bt$ is the optimal time of problem \thb{TI}   \thb{see \refeq{E109}}.
\end{Lemma}
\proof Let
$$
v_k(t)=b(t,u_k(t)), \qq t\in [0,T].
$$
By (P4), $v_k(\cd)$ is uniformly bounded in $L^\infty([0,T];\IR^n)$, which indicates that it is uniformly bounded in $L^2([0,T];\IR^n)$. Then $v_k(\cd)$ has sub-sequence, which converges weakly to some $v(\cd)$ in $L^2([0,T];\IR^n)$. Without loss of generality, let $v_k(\cd)$ itself weakly converges to $v(\cd)$ in $L^2([0,T];\IR^n)$.

Based on Mazur's Theorem  (see \cite{Yosida}, for example), there exists a sequence defined by the convex combination $\ds \sum^{N_k}_{j=1}\Ga_{k,j}v_j(\cd)$ strongly convergent to  $v(\cd)$ in $L^2([0,T];\IR^n)$. Since $U(t)$ is a compact convex set, we get
$$
v(t)\in U(t), \qq t\in [0,T].
$$
Then according to Filippov's Lemma (see \cite{Filippov} , for example), there exists $u(\cd)\in \mcU$, such that
$$
v(t)=b(t,u(t)),\qq \eqae t\in [0,T].
$$

Now, suppose $y(\cd)$ is the solution of equation \refeq{E106} with control $u(\cd)$. If $y(\cd)$ blows up during $[0,T]$, the lemma is proved. Otherwise,
$y(\cd)$ exists on $[0,T]$. Then, Lemma \ref{T201} shows that there exist $\Gd>0$ and $K>0$, such that $y_k(\cd)$ exists and is bounded uniformly on $[0,T+\Gd]$ when $k\geq K$, which contradicts to \refeq{E213}.

Therefore, $y(\cd)$ has to blow up in $[0,T]$, which proves our conclusion.
\endpf

Next, we can get the following existence theorem of problem (TI).
\begin{Theorem}\label{T202a}
Assume \thb{P1}---\thb{P4} hold and $\mcP_{ad}\ne \emptyset$. Then Problem \thb{TI} has at least one solution.
\end{Theorem}
\proof
Let $(T_k,y_k(\cd),u_k(\cd))\in\mcP_{ad}$ be a minimizing sequence, that is,
\begin{equation}\label{E214}
\lim_{k\to +\infty}T_k=\bt.
\end{equation}
Then
$$
T_k\geq \bt, \qq\all k\geq 1.
$$
Similar to the proof of Lemma \ref{T202}, there is $u(\cd)\in \mcU$, such that
$$
b(\cd,u_k(\cd))\to b(\cd,u(\cd)), \qq \mbox{weakly in}\,  L^2([0,\bt\,];\IR^n).
$$
Let $y(\cd)=y(\cd;u(\cd))$. Then we can easily see that $(\bt,y(\cd), u(\cd))\in \mcP$.

We claim that $y(\cd)$ blows up at $\bt$ \footnote{Based on the definition of $\bt$, $y(\cd;u(\cd))$ cannot blow up before $\bt$ .}, that is $(\bt,y(\cd), u(\cd))\in \mcP_{ad}$. Otherwise, by Lemma \ref{T201}, there exist $\Gd>0$ and $K>0$, such that $y_k(\cd)$ exists and is bounded uniformly on $[0,\bt+\Gd]$ for $k\geq K$. This contradicts \refeq{E214} and the fact that $y_k(\cd)$ blows up at $T_k$.

Therefore, $(\bt,y(\cd), u(\cd))\in \mcP_{ad}$, and $(\bt,y(\cd), u(\cd))$ is an optimal triple to Problem (TI).
\endpf

For Problem (TS), it holds that:
\begin{Theorem}\label{T203}
Assume \thb{P1}---\thb{P5} hold and $\mcU_{ad}=\mcU$. Let $t^*$ be defined by \thb{\ref{E110}} and it is finite. Then, Problem \thb{TS} has at least one solution.
\end{Theorem}
\proof
Let $(S_k,y_k(\cd),u_k(\cd))\in\mcP_{ad}$ be a maximizing sequence, that is
\begin{equation}\label{E215}
\lim_{k\to +\infty}S_k=t^*.
\end{equation}
Then
$$
S_k\leq t^*, \qq\all k\geq 1.
$$
Similar to the proof of Lemma \ref{T202}, there is $u(\cd)\in \mcU$, such that
$$
b(\cd,u_k(\cd))\to b(\cd,u(\cd)), \qq \mbox{weakly in}\,  L^2([0,t^*];\IR^n).
$$
Let $y(\cd)=y(\cd;u(\cd))$.

We claim that $y(\cd)$ blows up at $t^*$. Otherwise $y(\cd)$ blows up at some $S<t^*$ since $\mcU_{ad}=\mcU$.  By \refeq{E205}, there exists an $R>R_0$, such that
\begin{equation}\label{E216}
\int^\infty_R {1\over \zeta(r)}\, dr\leq {t^*-S\over 2}.
\end{equation}
Since $y(\cd)$ blows up at $S$, we have some $T<S$, such that
\begin{equation}\label{E217}
|y(T)|\geq  R+1.
\end{equation}
Using Lemma \ref{T201} and $\ds \lim_{k\to +\infty}S_k=t^*>S$, we get some $K>0$, such that $S_k>S$, and
\begin{equation}\label{E218}
|y_k(T)|\geq  R, \qq\all k\geq K.
\end{equation}
Noting that
\begin{eqnarray}\label{E219}
\nnb \doo{} t |y_k(t)| &=& G(t,|y_k(t)|)+{1\over |y_k(t)|}\ip{A(t)y_k(t)+b(t,u_k(t)),y_k(t)}\\
\nnb &\geq &  G(t,|y_k(t)|)-\|A(t)\| \, |y_k(t)|-\max_{u\in U} |b(t,u)|\\
&\geq & \zeta(|y_k(t)|), \qq t\in [T, S_k), \, k\geq K,
\end{eqnarray}
we have
\begin{eqnarray}\label{E220}
\nnb && {t^*-S\over 2}\geq  \int^{+\infty}_R{1\over \zeta(r)}\, dr\geq \int^{+\infty}_{|y_k(T)|}{1\over \zeta(r)}\, dr\\
&=& \int^{S_k}_T{1\over \zeta(|y_k(t)|)}\doo{} t |y_k(t)|\, dt \geq  S_k-T\geq S_k-S.
\end{eqnarray}
Let $k\to +\infty$, it follows that
\begin{equation}\label{E221}
{t^*-S\over 2}\geq t^*-S,
\end{equation}
which contradicts the assumption $t^*>S$. Thus, the blowup time of $y(\cd)$ is $t^*$. Therefore,
$(t^*,y(\cd), u(\cd))\in \mcP_{ad}$  and $(t^*,y(\cd), u(\cd))$ is an optimal triple of Problem (TS).
\endpf

\def\theequation{3.\arabic{equation}}
\setcounter{equation}{0} 
\setcounter{Definition}{0} \setcounter{Remark}{0}\section{Maximum Principles to Problem (TI)}

In this section, we will discuss the maximum principle of problem (TI). For simplicity, we may relabel some previous assumptions.

(S1) Let $(U,\rho)$ being a separable metric  space;

(S2) Function $G(t,r)$ is measurable in $t\in [0,+\infty)$, continuously differentiable in $r\in [0,+\infty)$ and
\begin{equation}\label{E301}
G(t,0)=0,   \qq  \all  t \in [0,+\infty).
\end{equation}
Moreover, $\all M>\Ga>0$,
\begin{equation}\label{E302}
\sup_{(t,r)\in [0,M]^2}|{\pa G(t,r)\over \pa r}|<+\infty,
\end{equation}
\begin{equation}\label{E303}
 \lim_{r\to +\infty}\inf_{t\in [\Ga,M]}{G(t,r)\over r}=+\infty,
\end{equation}
\begin{equation}\label{E304}
 \liminf_{r\to +\infty}\inf_{t\in [\Ga,M]}{rG_r(t,r)\over G(t,r)}>0.
\end{equation}

(S3) There exists $s_0>0$, $\Gvp(\cd) \in C^2(0,s_0)$ and modulus of continuity $\Go(\cd)\in C[0,+\infty)$, satisfying
\begin{equation}\label{E305}
\Gvp(s)>2, \q  \Gvp^\pri(s) <0, \qq\all s\in (0,s_0),
\end{equation}
\begin{equation}\label{E306}
 \lim_{s\to 0^+}\Gvp(s)=+\infty,\q \lim_{s\to 0^+}\Gvp^\pri(s)=-\infty,\q \lim_{s\to 0^+}{\Gvp(s)\over \Gvp^\pri(s)}=0,
\end{equation}
\begin{eqnarray}\label{E307}
  \nnb && 1+\Big|{\Gvp^\pri(s)\over \Gvp^2(s)}\Big|+  \Big|{\Gvp(s)\Gvp^\prii(s)\over (\Gvp^\pri(s))^2}\Big|\leq \Go(s)\Big(G_r(t,\Gvp(s))-{G(t,\Gvp(s))\Gvp^\prii(s)\over  (\Gvp^\pri(s))^2}\Big), \\
  && \qq \hspace{5cm} \all (t,s)\in [0,+\infty)\times (0,s_0).
\end{eqnarray}

(S4)  For any $T>0$,
\begin{equation}\label{E308}
 \sup_{t\in [0,T]} \|A(t)\|<+\infty, \q \sup_{t\in [0,T]\atop u\in U} |b(t,u)|<+\infty.
\end{equation}

To simplify the discussion, for$\rho>0$ and $s\in (0,s_0)$, denote
\begin{equation}\label{E309}
 \GO_T(\rho)=\inf_{t\in [0,T]\atop r\geq \rho} {G(t,r)\over r}, \, \Go_0(s)=\sup_{0<\tis<s}{1\over \Gvp(\tis)}, \, \Go_1(s)=\sup_{0<\tis<s}{\Gvp(\tis)\over |\Gvp^\pri(\tis)|},
\end{equation}
\begin{eqnarray*}
\tiGO(\rho)&=&\inf_{r\geq \rho}\Big(G_r(t,r)-{G(t,r)\Gvp^\prii(\Phi(r))\over  (\Gvp^\pri(\Phi(r)))^2}\Big)\\
&=&\inf_{\Gvp(s)\geq \rho}\Big(G_r(t,\Gvp(s))-{G(t,\Gvp(s))\Gvp^\prii(s)\over  (\Gvp^\pri(s))^2}\Big) ,
\end{eqnarray*}
where $\Phi(\cd)$ is the inverse function of $\Gvp(\cd)$.

For $f=\pmatrix{f^1 & f^2 & \ldots & f^n}^\top$, denote
$$
f_t=\pmatrix{{\pa f^1\over \pa t} \cr {\pa f^2\over \pa t}\cr \vdots\cr {\pa f^n\over \pa t}}, \q f_y= \dpp f y= \pmatrix{
  {\pa f^1\over \pa y_1} &  {\pa f^2\over \pa y_1} & \cdots & {\pa f^n\over \pa y_1}\cr
   {\pa f^1\over \pa y_2} &  {\pa f^2\over \pa y_2} & \cdots & {\pa f^n\over \pa y_2}\cr
   \vdots & \vdots & \ddots & \vdots \cr
     {\pa f^1\over \pa y_n} &  {\pa f^2\over \pa y_n} & \cdots & {\pa f^n\over \pa y_n}}.
$$


We have the following lemma.
\begin{Lemma}\label{T301} Assume \thb{S2}---\thb{S4} hold,  $T>t_0\geq 0$ and $g(\cd)\in L^\infty([t_0,T]; \IR^n)$. Let $\tiy(\cd)$ and $\hy(\cd)$ be the
solution of
\begin{equation}\label{E310}
 \doo {y(t)} t =G(t,|y(t)|){y(t)\over |y(t)|} +A(t)y(t)+g(t),\qq 　t>t_0
\end{equation}
on $[t_0,T]$ with the initial state $y(t_0)=\tiy_0$ and $y(t_0)=\hy_0$, respectively.

Suppose that
\begin{equation}\label{E309}
|\hy_0|>\rho,
\end{equation}
\begin{equation}\label{E311}
\Phi(|\hy_0|)-\Phi(|\tiy_0|)-\Big|{\hy_0\over |\hy_0|}-{\tiy_0\over |\tiy_0|}\Big|>0,
\end{equation}
where $\rho>0$  satisfies
\begin{equation}\label{E312}
\GO_T(\rho)\geq M+1,  \, \rho>2M, \, M\equiv\esssup_{t\in [t_0,T]}\max(|g(t)|,\|A(t)\|),
\end{equation}
$$
\Go_0(\Phi(\rho))\leq 1, \q \Go_1(\Phi(\rho))\leq 1,\q \Go(\Phi(\rho))\leq {1\over 4M+1},
$$
$$
 \tiGO(\rho)\geq 18M.
$$
Then, $\ds \Phi(|\hy(t)|)-\Phi(|\tiy(t)|)-\Big|{\hy(t)\over |\hy(t)|}-{\tiy(t)\over |\tiy(t)|}\Big|$ is monotonically increasing on $[t_0,T]$.
\end{Lemma}
\proof Since $\Phi(\cd)$ is monotonically decreasing, we get $|\tiy_0|>|\hy_0|>\rho$ from \refeq{E311}.

As a solution of \refeq{E310}, $y(\cd)$  satisfies
\begin{eqnarray}\label{E313}
\nnb \doo {|y(t)|} t&=& G(t,|y(t)|)+\ip{A(t)y(t)+g(t),{y(t)\over |y(t)|}}\\
&\geq & \Big(\GO_T(|y(t)|)-M\Big) |y(t)| - M, \qq\all t\in [t_0,T].
\end{eqnarray}
Hence, $|y(\cd)|$ is monotonically increasing on $[t_0,T]$  when $|y(t_0)|>\rho $. Especially,
\begin{equation}\label{E314}
|\tiy(t)|> \rho, \q |\hy(t)|>\rho, \qq\all t\in [t_0,T].
\end{equation}
Denote
$$
X(t)=|\hx(t)|-|\tix(t)|, \q \GT(t)= \hGt(t)-\tiGt(t),
$$
where
$$
\hGt(t)={\hy(t)\over |\hy(t)|},\q\tiGt={\tiy(t)\over |\tiy(t)|}, \q \hy(t)=\Gvp(|\hx(t)|) \hGt(t), \q \tiy(t)=\Gvp(|\tix(t)|) \tiGt(t),
$$
or equivalently,
$$
\hGt(t)={\hx(t)\over |\hx(t)|},\q\tiGt={\tix(t)\over |\tix(t)|},\q \hx(t)=\Phi(|\hy(t)|) \hGt(t), \q \tix(t)=\Phi(|\tiy(t)|) \tiGt(t).
$$
We have
\begin{eqnarray}\label{E315}
\nnb\ds \doo {X(t)} t &=& {G(t,\Gvp(|\hx(t)|))\over \Gvp^\pri(|\hx(t)|)}-{G(t,\Gvp(|\tix(t)|))\over \Gvp^\pri(|\tix(t)|)}+\ip{g(t),{\hGt(t)\over \Gvp^\pri(|\hx(t)|)} -{\tiGt(t)\over \Gvp^\pri(|\tix(t)|)} }   \\
\nnb && +{\Gvp(|\hx(t)|)\over \Gvp^\pri(|\hx(t)|)}\ip{A(t)\hGt(t),\hGt(t)} -{\Gvp(|\tix(t)|)\over \Gvp^\pri(|\tix(t)|)}\ip{A(t)\tiGt(t),\tiGt(t)}  \\
\nnb &=& {G(t,\Gvp(|\hx(t)|))\over \Gvp^\pri(|\hx(t)|)}-{G(t,\Gvp(|\tix(t)|))\over \Gvp^\pri(|\tix(t)|)}+\Big({1\over \Gvp^\pri(|\hx(t)|)} -{1\over \Gvp^\pri(|\tix(t)|)}\Big)\ip{g(t), \hGt(t)}\\
\nnb &&　+{1 \over \Gvp^\pri(|\tix(t)|)} \ip{g(t), \GT(t)}+\Big({\Gvp(|\hx(t)|)\over \Gvp^\pri(|\hx(t)|)} -{\Gvp(|\tix(t)|)\over \Gvp^\pri(|\tix(t)|)}\Big)\ip{A(t)\hGt(t), \hGt(t)}\\
\nnb &&　+{\Gvp(|\tix(t)|)\over \Gvp^\pri(|\tix(t)|)}\Big(\ip{A(t)\hGt(t),\GT(t)}+\ip{A(t)\GT(t),\tiGt(t)}\Big)\\
\nnb &\geq &  {G(t,\Gvp(|\hx(t)|))\over \Gvp^\pri(|\hx(t)|)}-{G(t,\Gvp(|\tix(t)|))\over \Gvp^\pri(|\tix(t)|)}\\
\nnb & & -M\Big|{1\over \Gvp^\pri(|\hx(t)|)} -{1\over \Gvp^\pri(|\tix(t)|)}\Big|- M\Big|{\Gvp(|\hx(t)|)\over \Gvp^\pri(|\hx(t)|)} -{\Gvp(|\tix(t)|)\over \Gvp^\pri(|\tix(t)|)}\Big|\\
 &&- 3M \Go_1(\Phi(\rho))\,|\GT(t)|,
\end{eqnarray}
\begin{eqnarray}\label{E316}
\nnb\ds \doo {\GT(t)} t 
\nnb &=& A(t)\GT(t)-\hGt(t)\hGt(t)^\top A(t)\hGt(t)+\tiGt(t)\tiGt(t)^\top A(t)\tiGt(t)\\
\nnb & &+{1\over \Gvp(|\hx(t)|)}\Big(g(t)- \hGt(t)\hGt(t)^\top g(t)\Big)-{1\over \Gvp(|\tix(t)|)}\Big(g(t)- \tiGt(t)\tiGt(t)^\top g(t)\Big)\\
\nnb &=& A(t)\GT(t)-\GT(t)\hGt(t)^\top A(t)\hGt(t)-\tiGt(t)\GT(t)^\top A(t)\hGt(t)-\tiGt(t)\tiGt(t)^\top A(t)\GT(t)\\
\nnb & &+\Big({1\over \Gvp(|\hx(t)|)}-{1\over \Gvp(|\tix(t)|)}\Big)\Big(g(t)- \hGt(t)\hGt(t)^\top g(t)\Big)\\
 && -{1\over \Gvp(|\tix(t)|)}\Big(\GT(t)\hGt(t)^\top g(t)+ \tiGt(t)\GT(t)^\top g(t)\Big).
\end{eqnarray}
Then,
\begin{eqnarray}\label{E317}
\ds \doo {|\GT(t)|} t &\leq &  M\Big(4+2\Go_0(\Phi(\rho)\Big)|\GT(t)|+2M\Big|{1\over \Gvp(|\hx(t)|)}-{1\over \Gvp(|\tix(t)|)}\Big|.
\end{eqnarray}
We know $X(t_0)-|\GT(t_0)|>0$. Denote
$$
S=\sup\set{\Gb\in (t_0,T]\Big| X(t)-|\GT(t)|>0, \q \all t\in [t_0,\Gb)}.
$$
Then,
$$
X(t)-|\GT(t)|>0, \qq \all t\in [t_0,S).
$$
Moreover, $X(S)-|\GT(S)|=0$ if $S<T$.

By \refeq{E315}---\refeq{E317}, we get
\begin{eqnarray}\label{E318}
\nnb \ds && \doo {\Big(X(t)-|\GT(t)|\Big)} t \\
\nnb &\geq &  {G(t,\Gvp(|\hx(t)|))\over \Gvp^\pri(|\hx(t)|)}-{G(t,\Gvp(|\tix(t)|))\over \Gvp^\pri(|\tix(t)|)}\\
\nnb && -M\Big|{1\over \Gvp^\pri(|\hx(t)|)} -{1\over \Gvp^\pri(|\tix(t)|)}\Big| -2M\Big|{1\over \Gvp(|\hx(t)|)}-{1\over \Gvp(|\tix(t)|)}\Big|　\\
\nnb && - M\Big|{\Gvp(|\hx(t)|)\over \Gvp^\pri(|\hx(t)|)} -{\Gvp(|\tix(t)|)\over \Gvp^\pri(|\tix(t)|)}\Big|\\
\nnb && -M\Big(4+2\Go_0(\Phi(\rho)) +3\Go_1(\Phi(\rho))\Big)\,|\GT(t)|\\
\nnb &=& \int^{|\hx(t)|}_{|\tix(t)|}\Big(G_r(t,\Gvp(s))-{G(t,\Gvp(s))\Gvp^\prii(s)\over (\Gvp^\pri(s))^2}\Big)\, ds\\
\nnb && - M\Big|\int^{|\hx(t)|}_{|\tix(t)|}\Big(1-{\Gvp(s)\Gvp^\prii(s)\over (\Gvp^\pri(s))^2}\Big) \, ds\Big|\\
\nnb && -M \Big|\int^{|\hx(t)|}_{|\tix(t)|}{\Gvp^\prii(s)\over (\Gvp^\pri(s))^2} \, ds\Big|-2M\Big|\int^{|\hx(t)|}_{|\tix(t)|}{\Gvp^\pri(s)\over (\Gvp(s))^2} \, ds\Big|\\
\nnb && -M\Big(4+2\Go_0(\Phi(\rho)) +3\Go_1(\Phi(\rho))\Big)\,|\GT(t)|\\
\nnb &\geq & \int^{|\hx(t)|}_{|\tix(t)|}(1-2M\Go(\Phi(\rho))\Big(G_r(t,\Gvp(s))-{G(t,\Gvp(s))\Gvp^\prii(s)\over (\Gvp^\pri(s))^2}\Big)\, ds -7M\,|\GT(t)|\\
\nnb &\geq & {1\over 2}\tiGO(\rho)X(t)-9M\,|\GT(t)|\\
\nnb &\geq & 9M\,(X(t)-|\GT(t)|)\\
&\geq & 0, \qq t\in [t_0,S] .
\end{eqnarray}
Therefore, $X(t)-|\GT(t)|$ is monotonically increasing on $[t_0,S]$. Consequently
 $X(S)-|\GT(S)|>0$. Thus, $S=T$. We get the proof.
\endpf

Thanks the above lemma, we can
 obtain the maximum principle of problem (TI) easily now.
\begin{Theorem}\label{T302}
Assume  that \thb{S1}---\thb{S4} hold. Let $(\bt,\by(\cd),\bu(\cd))$ be an optimal triple of \thb{TI}.
Then, there exists a nontrivial solution  $\bpsi(\cd)\in C([0,\bt\,);\IR^n)$ of the following equation
\begin{equation}\label{E319}
\doo {\bpsi(t)} t=-\Big({G(t,|\by(t)|)\over |\by(t)|}I+{|\by(t)|\, G_r(t,|\by(t)|)-G(t,|\by(t)|)\over |\by(t)|^3}\by(t)\by(t)^\top+A(t)^\top\Big)\bpsi(t),\q    t\in [0,\bt\,)
\end{equation}
such that
\begin{eqnarray}\label{E320}
 &&  \ip{\bpsi(t),   b(t,\bu(t))} =\max_{ u\in U}\ip{\bpsi(t),  b(t,u)},\q \eqae t\in [0,\bt\,).
\end{eqnarray}
Moreover,
\begin{equation}\label{E321}
   \lim_{t\to \bt^-}\bpsi(t)=0.
\end{equation}
\end{Theorem}
\proof Denote
$$
M=\esssup_{t\in [0,\bt\,]} \sup_{u\in U}\max(|b(t,u)|,\|A(t)\|).
$$
Set $\rho>0$ such that
$$
\GO_\bt(\rho)\geq M+1,  \q \rho>2M,
\q  \Go_0(\Phi(\rho))\leq 1,\q  \Go_1(\Phi(\rho))\leq 1,  \q \Go(\Phi(\rho))\leq {1\over 4M+1},
$$
$$
 \tiGO(\rho)\geq 18M.
$$
It is not difficult to see the existence of such $\rho$. On the other hand, there exists $\Gd>0$, such that
$$
|\by(t)|\geq \rho, \qq t\in [\bt-\Gd,\bt\,).
$$
For $|z|>\rho$, let
$$
E_z\equiv \set{\ell z|  \ell\geq 1}.
$$
Then, by Lemma \ref{T301}, $(\by(\cd),\bu(\cd))$ is an optimal pair of the following optimal control problem for any $T\in [\bt-\Gd,\bt\,)$: to find a control $u(\cd)\in \mcU$, such that the solution $y(\cd)$ of
\begin{equation}\label{E322}
 \left\{\begin{array}{ll}\ds  \doo {y(t)} t =G(t,|y(t)|){y(t)\over |y(t)|}+A(t)y(t)+b(t,u(t)),&　t\in [0,T],\vspace{4mm}\\
   y(0)=y_0  & \end{array}\right.
\end{equation}
maximizes $|y(T)|^2$ with terminal constraint $y(T)\in E_{\by(T)}$.

Otherwise, there exists $\tiu(\cd)\in \mcU$ and $\ell>1$, such that
$$
y(T;\tiu(\cd))=\ell \by(T).
$$
In this case,
$$
\Phi(|\by(T)|)-\Phi(|y(T;\tiu(\cd)|)- \Big|{\by(T)\over |\by(T)|}-{y(T;\tiu(\cd)\over |y(T;\tiu(\cd)|}\Big|=\Phi(|\by(T)|)-\Phi(\ell |\by(T)|)>0.
$$
Define
$$
\hu(t)= \left\{\begin{array}{ll}\ds  \tiu(t), & t\in [0,T], \hspace{4mm}\\
   \ds  \bu(t), & t\in [T,\bt\,). \end{array}\right.
$$
By Lemma \ref{T301}, $\Phi(|\by(\cd)|)-\Phi(|y(\cd;\hu(\cd)|)- \Big|{\by(\cd)\over |\by(\cd)|}-{y(\cd;\hu(\cd)\over |y(\cd;\hu(\cd)|}\Big|$ increases in the existence interval of $y(\cd;\hu(\cd))$ within  $[T, \bt)$. In particular, we have some $S<\bt$, satisfying $\ds \lim_{t\to S^-}\Phi(|y(t;\hu(\cd)|)=0$. That is
$y(\cd;\hu(\cd))$ blows up at $S$, which contradicts the optimality of $(\bt,\by(\cd),\bu(\cd))$.

Then, using the classical maximum principle, there exists a nontrivial pair $(\Gvp_{0,T}, \Gvp_T(\cd))\in \IR\times C([0,T];\IR^n)$, which satisfies
$$
\Gvp_{0,T}\leq 0,
$$
\begin{equation}\label{E323}
 \doo {\Gvp_T(t)} t=-\Big({G(t,|\by(t)|)\over |\by(t)|}I+{|\by(t)|\, G_r(t,|\by(t)|)-G(t,|\by(t)|)\over |\by(t)|^3}\by(t)\by(t)^\top+A(t)^\top\Big)\Gvp_T(t),\q   t\in [0,T],
\end{equation}
\begin{equation}\label{E324}
\ip{\Gvp_T(t), b(t,\bu(t))} =\max_{u\in U}\ip{\Gvp_T(t),b(t,u)},\q \eqae t\in [0,T]
\end{equation}
and
\begin{equation}\label{E325}
\ip{\Gvp_T(T)+\Gvp_{0,T}\by(T),q-\by(T)}\geq 0, \qq\all q\in E_{\by(T)}.
\end{equation}
Obviously, \refeq{E325} implies
\begin{equation}\label{E326}
\ip{\Gvp_T(T)+\Gvp_{0,T}\by(T), \by(T)}\geq 0.
\end{equation}
If $\Gvp_{0,T}=0$, $\Gvp_T(\cd)\ne 0$ because of the  non-triviality. If $\Gvp_{0,T}\ne 0$, then it follows from \refeq{E326} that
$$
\ip{\Gvp_T(t), \by(T)}\geq  -\Gvp_{0,T}|\by(T)|^2>0,
$$
which indicates $\Gvp_T(\cd)\ne 0$, either.

Therefore, $\Gvp_T(\cd)\ne 0$ always holds. Thus, replacing $\Gvp_T(\cd)$ by $\ds {\Gvp_T(\cd)\over |\Gvp_T(0)|}$ if necessary,  we
can suppose that $|\Gvp_T(0)|=1$.  While \refeq{E323}---\refeq{E326} remain true.

Now, for any $\Gve>0$, we can see that $\Gvp_T(\cd)$ is equicontinuous on $[0,\bt-\Gve]$. Therefore $\Gvp_T(\cd)$  has a subsequence that converges uniformly to  $\bpsi(\cd)$ on $[0,\bt-\Gve]$ when $T\to \bt^-$. Let $\Gvp_T(\cd)$ itself be the subsequence for simplicity. Then, we get \refeq{E319}---\refeq{E320}.

Furthermore, by \refeq{E304}, we know there exist $T_0\in (0,\bt\,)$ and $c\in (0,1)$  such that
$$
|\by(t)|\, G_r(t,|\by(t)|)\geq c \,G(t,|\by(t)|), \qq \all t\in [T_0,\bt\,).
$$
By \refeq{E319}, we get
\begin{eqnarray}\label{E327}
\nnb  {1\over 2}\doo{|\bpsi(t)|^2} t& =& - {G(t,|\by(t)|)\over |\by(t)|}|\bpsi(t)|^2-{|\by(t)|\, G_r(t,|\by(t)|)-G(t,|\by(t)|)\over |\by(t)|^3}\ip{\by(t),\bpsi(t)}^2\\
\nnb &&  -\ip{A(t)\bpsi(t), \bpsi(t)}\\
 &\leq & -\Big( {c\, G(t,|\by(t)|)\over |\by(t)|}-M\Big)\,|\bpsi(t)|^2,\q    t\in [T_0,\bt\,).
\end{eqnarray}
On the other hand, using
\begin{eqnarray}\label{E328}
\nnb  && {1\over 2}\doo{ |\by(t)|^2} t =  {G(t, |\by(t)|)\over |\by(t)|}\, |\by(t)|^2+\ip{A(t)\by(t), \by(t)}+{b(t,\bu(t)),\by(t)}\\
\nnb &\leq &   \Big({G(t, |\by(t)|)\over |\by(t)|}+M\Big)\, |\by(t)|^2 +M|\by(t)|, \qq\all t\in [0,\bt\,)
\end{eqnarray}
and $\ds \lim_{t\to \bt^-}|\by(t)|=+\infty$, we conclude
$$
\lim_{t\to \bt^-}\int^t_0 \Big({G(s, |\by(s)|)\over |\by(s)|}+M\Big)\, ds=+\infty.
$$
Therefore, by \refeq{E326}
$$
\lim_{t\to \bt^-}\int^t_{T_0} {G(s, |\by(s)|)\over |\by(s)|}\, ds=+\infty.
$$
Thus,  \refeq{E321} can be derived from \refeq{E327} and the above inequality. \endpf

We will find later that any optimal triple of problem (TS) also satisfies the above theorem. Thus, we would like to make some further observation on the optimal triple of problem (TI). We assume that:

(S5) Let $b(t,U)$ be a convex set with the origin point being its interior point for almost all $t\in [0,+\infty)$. Meanwhile, for any $x\in \pa \Big(b(t,U)\Big)$, there exists a unique  $\Gl\in S^{n-1}$, such that
\begin{equation}\label{E327}
\ip{\Gl, y-x}\leq 0, \qq \all y\in b(t,U).
\end{equation}

\begin{Remark}
If $U$ is a closed ball in $\IR^m$, in which the origin is an interior point,
\begin{equation}\label{E330}
b(t,u)=B(t)u, \qq \all t\in [0,+\infty), u\in U,
\end{equation}
and $B(t)\in \IR^{n\times m}$ always has full row rank, then \thb{S5} holds.
\end{Remark}

We have:
\begin{Theorem}\label{T303}
Assume that \thb{S1}---\thb{S5} hold  and $(\bt,\by(\cd),\bu(\cd))$ is an optimal triple of problem \thb{TI}.
Then, there exists a nontrivial solution $\bpsi(\cd)\in C([0,\bt\,);\IR^n)$ of equation \refeq{E319}, such that \refeq{E320}---\refeq{E321} hold. Moreover, we have some $\Gd\in (0,\bt\,)$ and the following transversality condition:
\begin{equation}\label{E331}
\ip{\bpsi(t), \by(t)}> 0, \qq \all t\in (\bt-\Gd,\bt\,).
\end{equation}
\end{Theorem}
\proof We will use symbols that used in the proof of Theorem \ref{T302}. We need only to prove \refeq{E331}. If $ \bt-\Gd\leq  T_1<T_2<\bt$, then it follows from \refeq{E324} that
\begin{equation}\label{E332}
\ip{\Gvp_{T_i}(t), b(t,u)-b(t,\bu(t))}\leq 0, \qq \eqae t\in [0,T_1], \q i=1,2.
\end{equation}
Since $\Gvp_{T_i}(t)\ne 0$ ($\all t\in [0,T_i]$, $i=1,2$), we find out $b(t,\bu(t))$ is a boundary point of $b(t,U)$ for almost all $t\in [0,T_i]$. Thus, (S5) and \refeq{E332} imply that
\begin{equation}\label{E333}
\Gvp_{T_1}(t)=c\Gvp_{T_2}(t), \qq \eqae t\in [0,T_1]
\end{equation}
for some constant $c>0$.
Using the continuity of $\Gvp_{T_i}(\cd)$ on $[0,T_1]$ and $|\Gvp_{T_i}(0)|=1$ ($i=1,2$), we get
\begin{equation}\label{E334}
\Gvp_{T_1}(t)=\Gvp_{T_2}(t), \qq  t\in [0,T_1].
\end{equation}
Consequently,
\begin{equation}\label{E335}
\bpsi (t)=\Gvp_T(t), \qq \all  t\in [0,T], \, T\in [\bt-\Gd,\bt\,)
\end{equation}
since $\Gvp_T(\cd)$ converges uniformly to $\bpsi(\cd)$ on $[0,\bt-\Gve]$ for any $\Gve>0$.
Especially,
\begin{equation}\label{E336}
\ip{\bpsi (T), \by(T)}=\ip{\Gvp_T(T),\by(T)} \geq -\Gvp_{0,T}|\by(T)|^2\geq 0, \qq \all T\in [\bt-\Gd,\bt\,).
\end{equation}
On the other hand, since $0$ is an interior point in $b(t,U)$ for almost all $t\in [0,\bt\,)$,
\begin{equation}\label{E337}
\ip{\bpsi (t), b(t,\bu(t))}=\max_{u\in U}\ip{\bpsi (t), b(t,u)}> 0, \qq \eqae t\in (0,\bt\,).
\end{equation}
Then
\begin{eqnarray}\label{E338}
\nnb   && \doo {}t  \ip{\bpsi (t), \by(t)} \\
\nnb &=& - {|\by(t)|\, G_r(t,|\by(t)|)-G(t,|\by(t)|)\over |\by(t)|}  \ip{\bpsi(t), \by(t)} + \ip{\bpsi (t),  b(t,\bu(t))}\\
&>&  - {|\by(t)|\, G_r(t,|\by(t)|)-G(t,|\by(t)|)\over |\by(t)|}  \ip{\bpsi(t), \by(t)}, \qq \eqae t \in (0,\bt\,).
\end{eqnarray}
Using (S2),
\begin{equation}\label{E339}
g(t)\equiv \int^t_0 {|\by(s)|\, G_r(s,|\by(s)|)-G(s,|\by(s)|)\over |\by(s)|} \, ds
\end{equation}
is well-defined in $(0,\bt\,)$ and it follows from \refeq{E338} that
\begin{equation}\label{E340}
\doo {}t \Big(e^{g(t)} \ip{\bpsi (t), \by(t)}\Big)>0, \qq\all t\in (0,\bt\,).
\end{equation}
Finally, \refeq{E331} follows easily from \refeq{E336} and \refeq{E340}. We complete the proof.
\endpf

\def\theequation{4.\arabic{equation}}
\setcounter{equation}{0} 
\setcounter{Definition}{0} \setcounter{Remark}{0}\section{Maximum Principles to Problem (TS)}

For Problem (TS), we assume that

(S5$^\prime$) For almost all $t\in [0,+\infty)$, the origin point is an interior point  of
$b(t,U)$.

We have:
\begin{Theorem}\label{T304}
Assume that \thb{S1}---\thb{S4} hold  and $(t^*,y^*(\cd),u^*(\cd))$ is an optimal triple of Problem \thb{TS}.
Then there exists a nontrivial solution $\psi^*(\cd)\in C([0,t^*);\IR^n)$ of the following equation
\begin{equation}\label{E401}
\doo {\psi^*(t)} t=-\Big({G(t,|y^*(t)|)\over |y^*(t)|}I+{|y^*(t)|\, G_r(t,|y^*(t)|)-G(t,|y^*(t)|)\over |y^*(t)|^3}y^*(t)y^*(t)^\top+A(t)^\top\Big)\psi^*(t),\q    t\in [0,t^*)
\end{equation}
such that
\begin{eqnarray}\label{E402}
 &&  \ip{\psi^*(t),   b(t,u^*(t))} =\max_{ u\in U}\ip{\psi^*(t),  b(t,u)},\q \eqae t\in [0,t^*)
\end{eqnarray}
and
\begin{equation}\label{E403}
   \lim_{t\to t^{*-}}\psi^*(t)=0.
\end{equation}
Furthermore, if \thb{S5$^\prime$} holds, then
\begin{equation}\label{E404}
\ip{\psi^*(t), y^*(t)}<0, \qq \all t\in [0,t^*).
\end{equation}
\end{Theorem}
\proof The proof is similar to that of Theorems \ref{T302} and \ref{T303}.
Let
$$
M=\esssup_{t\in [0,t^*]} \sup_{u\in U}\max(|b(t,u)|,\|A(t)\|)
$$
and $\rho>0$, such that
$$
\GO_{t^*}(\rho)\geq M+1,  \q \rho>2M,
\q  \Go_0(\Phi(\rho))\leq 1,\q  \Go_1(\Phi(\rho))\leq 1,  \q \Go(\Phi(\rho))\leq {1\over 4M+1},
$$
$$
 \tiGO(\rho)\geq 18M.
$$
The existence of such a $\rho$ is obvious, and we have some $\Gd>0$, such that
$$
|y^*(t)|\geq 2\rho, \qq \all t\in [t^*-\Gd,t^*).
$$
For $|z|>\rho$,  denote\footnote{\zihao{6} In the proof of Theorem \ref{T302}, it will also work if $E_z$ was defined by \refeq{E405} there.}
\begin{equation}\label{E405}
E_z\equiv \set{\ell z\Big| {1\over 2}\leq \ell\leq 1}.
\end{equation}
Then, for any $T\in [t^*-\Gd,t^*)$, it is easy to derive that $(y^*(\cd),u^*(\cd))$ is an optimal pair of the following problem by Lemma \ref{T301}: to find a control $u(\cd)\in \mcU$, such that the solution $y(\cd)$ of
\begin{equation}\label{E406}
 \left\{\begin{array}{ll}\ds  \doo {y(t)} t =G(t,|y(t)|){y(t)\over |y^*(t)|}+A(t)y(t)+b(t,u(t)),&　t\in [0,T],\vspace{4mm}\\
   y(0)=y_0  & \end{array}\right.
\end{equation}
minimizes $|y(T)|^2$ with the terminal constraint $y(T)\in E_{y^*(T)}$.

Otherwise, there exists a $\tiu(\cd)\in \mcU$ and an $\ell\in [{1\over 2}, 1)$, such that
$$
y(T;\tiu(\cd)=\ell y^*(T).
$$
In this case
$$
\Phi(|y(T;\tiu(\cd)|)-\Phi(|y^*(T)|)-\Big|{y(T;\tiu(\cd)\over |y(T;\tiu(\cd)|}-{y^*(T)\over |y^*(T)|}\Big|=\Phi(\ell |y^*(T)|)-\Phi(|y^*(T)|)>0.
$$
We set
$$
\hu(\cd)= \left\{\begin{array}{ll}\ds  \tiu(t), & t\in [0,T], \hspace{4mm}\\
   \ds  u^*(t), & t\in [T,t^*). \end{array}\right.
$$
Then, by Lemma \ref{T301}, $\Phi(|y(\cd;\hu(\cd)|)-\Phi(|y^*(\cd)|)- \Big|{y(\cd;\hu(\cd)\over |y(\cd;\hu(\cd)|}-{y^*(\cd)\over |y^*(\cd)|}\Big|$ is monotonically increasing in the existence interval of $y(\cd;\hu(\cd))$ within $[T, t^*)$. In particular, we know that $y(\cd;\hu(\cd))$
exists on $[0,t^*]$ and $\ds\lim_{t\to t^{*-}}\Phi(|y(t^*;\hu(\cd)|)>0$, which contradicts the optimality of $(t^*,y^*(\cd),u^*(\cd))$.

According to the classical maximum principle, we get a nontrivial pair $(\Gvp_{0,T}, \Gvp_T(\cd))\in \IR\times C([0,T];\IR^n)$, satisfying
$$
\Gvp_{0,T}\leq 0,
$$
\begin{equation}\label{E407}
 \doo {\Gvp_T(t)} t=-\Big({G(t,|y^*(t)|)\over |y^*(t)|}I+{|y^*(t)|\, G_r(t,|y^*(t)|)-G(t,|y^*(t)|)\over |y^*(t)|^3}y^*(t)y^*(t)^\top+A(t)^\top\Big)\Gvp_T(t),\q   t\in [0,T],
\end{equation}
\begin{equation}\label{E408}
\ip{\Gvp_T(t), b(t,u^*(t))} =\max_{u\in U}\ip{\Gvp_T(t),b(t,u)},\q \eqae t\in [0,T].
\end{equation}
and
\begin{equation}\label{E409}
\ip{\Gvp_T(T)-\Gvp_{0,T}y^*(T),q-y^*(T)}\geq 0, \qq\all q\in E_{y^*(T)}.
\end{equation}
Obviously, \refeq{E409} ensures
\begin{equation}\label{E410}
\ip{\Gvp_T(T), y^*(T)}= \Gvp_{0,T}|y^*(T)|^2.
\end{equation}
If $\Gvp_{0,T}=0$, we get $\Gvp_T(\cd)\ne 0$ from the non-triviality. If $\Gvp_{0,T}\ne 0$,  we get from \refeq{E410} that
$$
\ip{\Gvp_T(t), y^*(T)}<0.
$$
Then, $\Gvp_T(\cd)\ne 0$ also holds.

To conclude, $\Gvp_T(\cd)\ne 0$ is always tenable. Then, we can reset  $\Gvp_T(\cd)$ such that $|\Gvp_T(0)|=1$. While  \refeq{E407}---\refeq{E408} and
\begin{equation}\label{E411}
\ip{\Gvp_T(T), y^*(T)}\leq 0
\end{equation}
still hold.

Next, similar to the proof of Theorem \ref{T302}, we get that,  at least along a subsequence, $\Gvp_T(\cd)$ convergence uniformly to
 $\psi^*(\cd)$ on $[0,t^*-\Gve]$ for any $\Gve>0$. Then we get the conjugate function $\psi^*(\cd)$ and \refeq{E401}---\refeq{E403}.

When (S5$^\prime$) holds, similar to \refeq{E340},  we have
\begin{equation}\label{E412}
\doo {}t \Big(e^{h(t)} \ip{\Gvp_T (t), \by(t)}\Big)>0, \qq\eqae t\in (0,T), \, T\in [t^*-\Gd,t^*),
\end{equation}
where
\begin{equation}\label{E413}
h(t)\equiv \int^t_0 {|y^*(s)|\, G_r(s,|y^*(s)|)-G(s,|y^*(s)|)\over |y^*(s)|} \, ds, \qq t\in [0,t^*).
\end{equation}
Combining\refeq{E411} with \refeq{E412}, we get
\begin{equation}\label{E414}
\ip{\Gvp_T(t), y^*(t)}< 0, \qq\eqae t\in [0,T).
\end{equation}
Therefore,
\begin{equation}\label{E415}
\ip{\psi^*(t), y^*(t)}\leq 0, \qq\all t\in [0,t^*).
\end{equation}
Then, since it also holds that
\begin{equation}\label{E416}
\doo {}t \Big(e^{h(t)} \ip{\psi^*(t), \by(t)}\Big)>0, \qq\eqae t\in [0,t^*),
\end{equation}
we get \refeq{E404}.
\endpf

\begin{Remark}\label{R304} In the proof of Theorem \ref{T304}, if \thb{S5$^\pri$} is replaced by the following condition:

$(S5^\prii)$  For almost all $t\in [0,+\infty)$,
\begin{equation}\label{E417}
0\in \overline{b(t,U)},
\end{equation}
then in stead of \refeq{E404},  we would get the following transversality condition:
\begin{equation}\label{E418}
\ip{\psi^*(t), y^*(t)}\leq 0, \qq \all t\in [0,t^*).
\end{equation}
\end{Remark}

\def\theequation{5.\arabic{equation}}
\setcounter{equation}{0} 
\setcounter{Definition}{0} \setcounter{Remark}{0}\section{Results for Autonomous Systems.}

If \refeq{E106} is an autonomous system, then maximum principles and their proofs can be simplified. That is because in such a situation, when we limit the optimal triple $(\bt,\by(\cd),\bu(\cd))$ of problem(TI)/(TS) on $[0,T]$ for any $T\in (0,\bt\,)$, it should be a solution of a time optimal control problem that changes the state from $y_0$ to the target set $\set{\by(T)}$ most quickly/slowly. In this case, we can conclude the maximum principle of Problem (TI)/(TS) from the classical results on $[0,T]$ by taking $T\to \bt^-$.

Specifically, consider
\begin{equation}\label{E501}
   \left\{\begin{array}{ll}\ds  \doo {y(t)} t=f(y(t),u(t)),&　t>0,\\
   y(0)=y_0. & \end{array}\right.
\end{equation}
We set the following assumptions:

(A1) Let $(U,\rho)$ being a separable metric  space;

(A2) Function $f(y,u)$ is continuous in $(y,u)$ and continuously differentiable in $y \in \IR^n$. Meanwhile,
\begin{equation}\label{E502}
   |y(0,u)| \leq L, \qq\all u \in U
\end{equation}
for some $L>0$. Moreover, for any $R>0$, there exists $L_R>0$ such that
\begin{equation}\label{E503}
 | f_y(y,u)| \leq L_R,\qq \all |y|\leq R; u \in U.
\end{equation}

\begin{Theorem}\label{T501}
 Assume that \thb{A1}--\thb{A2} hold and $(\bt, \by(\cd), \bu(\cd))$ is an optimal triple of Problem \thb{TI}/\thb{TS}. Then, there exists a nontrivial solution  $\bpsi(\cd)\in C([0,\bt\,);\IR^n)$ of the following equation
\begin{equation}\label{E504}
\doo {\bpsi(t)} t=-f_y(\by(t), \bu(t))\bpsi(t),\q    t\in [0,\bt\,),
\end{equation}
such that
\begin{equation}\label{E505}
\ip{\bpsi(t), f(\by(t), \bu(t))} =\max_{ u\in U}\ip{\bpsi(t), f(\by(t), u)},\q \eqae t\in [0,\bt\,),
\end{equation}
\end{Theorem}
\proof The translation invariance of the autonomous systems ensures that for any $ T \in (0, \bt\,)$,  $(T, \by(\cd), \bu(\cd))$ is an optimal triple of the following optimal control problem:
to find  $(t^*,y^*(\cd),u^*(\cd))\in \mcP^T_{ad}$  such that
\begin{equation}\label{E506}
t^*=\inf_{(t,y(\cd),u(\cd))\in \mcP^T_{ad}}t \q\Big/\q  t^*=\sup_{(t,y(\cd),u(\cd))\in \mcP^T_{ad}}t,
\end{equation}
where
\begin{equation}\label{E507}
   \mcP^T_{ad}=\set{(t,y(\cd),u(\cd))\in (0,t]\times C([0,+\infty);\IR^n)\times \mcU\Big| \mbox{\refeq{E501} \,holds on }\, [0,t), \, y(t)=\by(T)}.
\end{equation}
Thus, by classical maximum principle, there exists a nontrivial solution  $\bpsi_T(\cd)\in C([0,T];\IR^n)$ of the following equation
\begin{equation}\label{E508}
\doo {\bpsi_T(t)} t=-f_y(\by(t), \bu(t))\bpsi_T(t),\qq    t\in [0,T]
\end{equation}
such that
\begin{eqnarray}\label{E509}
 &&  \ip{\bpsi_T(t), f(\by(t), \bu(t))} =\max_{ u\in U}\ip{\bpsi_T(t), f(\by(t), u)},\q \eqae t\in [0,T].
\end{eqnarray}

Because of the non-triviality of $\bpsi_T(\cd)$, we  can set $|\bpsi_T(0)|=1$. Thus, for any $\varepsilon >0$, $\bpsi_T(\cd)$ is equicontinuous on $[0,\bt-\varepsilon]$. So when $T\to \bt^-$, $\Gvp_T(\cd)$  has a subsequence that converges uniformly to  $\bpsi(\cd)$ on $[0,\bt-\Gve]$ for any $\Gve>0$. We get \refeq{E504}---\refeq{E505} and complete the  proof.
\endpf

\begin{Remark}\label{R501} In the assumptions of Theorem \ref{T501}, we do not set any nonlinear growth condition. This does not mean that such conditions are not important. Instead, these conditions are implied in the existence of an optimal triple.
\end{Remark}
\begin{Remark}\label{R502} We need some mild assumptions to yield transversality conditions. We will not discuss transversality conditions since it is quite  technical when assumptions are weak.
\end{Remark}

\def\theequation{6.\arabic{equation}}
\setcounter{equation}{0} 
\setcounter{Definition}{0} \setcounter{Remark}{0}\section{Some Examples}

We mention that the assumptions in Theorems \ref{T303} and \ref{T304} concerns mainly about $G$ and $\Gvp$. They looks quite technical. In fact, most systems that we care about satisfy these assumptions---(S2) and (S3). We list them in the following.

We always assume that $p>1$, $\ds \Gb>{1\over p-1}$ and $g(\cd)$ is measurable in $[0,+\infty)$, satisfying
\begin{equation}\label{E601}
0<\essinf_{t\in [\Ga,T]} g(t)\leq \esssup_{t\in [\Ga,T]} g(t)<+\infty, \qq\all T>\Ga>0.
\end{equation}

\begin{Example} Let
$$
G(t,r)=g(t) r^p, \qq t\in [0,+\infty), r\geq 0.
$$
In this case, take
\begin{equation}\label{E602}
\Gvp(s)=s^{-\Gb}, \qq s>0.
\end{equation}
It can be verified directly that \thb{S2}---\thb{S3} hold.
\end{Example}

\begin{Example} Let
$$
G(t,r)=g(t) r\ln^p(1+r), \qq t\in [0,+\infty), r\geq 0.
$$
We take
\begin{equation}\label{E603}
\Gvp(s)=\exp \big(s^{-\Gb}\big),   \qq s>0.
\end{equation}
Then \thb{S2}---\thb{S3} hold.
\end{Example}
\begin{Example} Let
$$
G(t,r)=g(t) e^{(p-1)r}, \qq t\in [0,+\infty), r\geq 0.
$$
Choose
\begin{equation}\label{E604}
\Gvp(s)=\ln (1+s^{-\Gb}), \qq s>0.
\end{equation}
Then \thb{S2}---\thb{S3} hold.
\end{Example}

\begin{Remark} In Examples 1 and 3, $\Gvp(\cd)$ can be defined by \refeq{E603},  one need only the positivity of $\Gb$ to guarantee \thb{S2}---\thb{S3}.

Similarly, in Example 3, if $\Gb>0$ and $\Gvp(\cd)$ is defined by \refeq{E602}, then \thb{S2}---\thb{S3} hold.
\end{Remark}

\end{document}